\documentclass[reqno]{amsart}

\usepackage{amsmath,amssymb,amsthm,mathtools}
\usepackage{enumitem}
\usepackage[a4paper,margin=32mm]{geometry}
\usepackage{microtype,url}
\usepackage[colorlinks=true,linkcolor=blue,citecolor=blue,urlcolor=blue]{hyperref}

\title[Parallelepipeds in ellipsoids]{Parallelepipeds of maximal facet area and total edge length in ellipsoids, through prescribed boundary points}
\author[T.~Kania]{Tomasz Kania}
\address[T.~Kania]{Mathematical Institute\\ Czech Academy of Sciences\\ \v Zitn\'a 25\\ 115 67 Praha 1\\ Czech Republic\\
Institute of Mathematics and Computer Science\\ Jagiellonian University\\ {\L}ojasiewicza 6\\ 30-348 Krak\'ow, Poland}
\email{tomasz.marcin.kania@gmail.com}
\thanks{IM CAS (RVO 67985840). Corresponding author: Tomasz Kania; \url{tomasz.marcin.kania@gmail.com}}
\keywords{Ellipsoid, inscribed parallelepiped, facet area, total edge length, Schur--Horn theorem, Connes--Zagier property}
\subjclass[2020]{Primary 52A40; Secondary 52B11, 52A38, 15A42, 15A45}
\date{\today}

\numberwithin{equation}{section}
\theoremstyle{plain}
\newtheorem{theorem}{Theorem}[section]
\newtheorem{proposition}[theorem]{Proposition}
\newtheorem{lemma}[theorem]{Lemma}
\newtheorem{corollary}[theorem]{Corollary}
\theoremstyle{definition}
\newtheorem{definition}[theorem]{Definition}
\theoremstyle{remark}
\newtheorem{remark}[theorem]{Remark}
\newtheorem{problem}[theorem]{Problem}

\newcommand{\R}{\mathbb{R}}
\newcommand{\tr}{\operatorname{tr}}
\newcommand{\diag}{\operatorname{diag}}
\newcommand{\Diag}{\operatorname{Diag}}
\newcommand{\Id}{\mathrm{I}}
\newcommand{\cE}{\mathcal{E}}
\newcommand{\one}{\mathbf{1}}
\newcommand{\ip}[2]{\langle #1,#2\rangle}

\renewcommand{\le}{\leqslant}
\renewcommand{\ge}{\geqslant}

\begin{document}

\begin{abstract}
Let
\[
   \cE_A=\{x\in\R^n:x^{\top}A^{-1}x\le 1\},\qquad n\ge2,
\]
where \(A\) is real symmetric positive definite.  We study full-dimensional
parallelepipeds whose \(2^n\) vertices lie on \(\partial\cE_A\).  First we show that
such parallelepipeds are necessarily centred at the origin and are precisely the
images, under \(A^{1/2}\), of orthotopes inscribed in the Euclidean unit sphere.
This reduces the extremal questions to finite-dimensional linear algebra.

For the total length \(L\) of the one-skeleton we prove
\[
       L_{\max}(\cE_A)=2^n\sqrt{\tr A}.
\]
Moreover, the prescribed-vertex problem for \(L\) has the same answer in every
dimension: for every \(x_0\in\partial\cE_A\) there is an inscribed parallelepiped
with vertex \(x_0\) and total edge length \(2^n\sqrt{\tr A}\).  The proof uses the
Schur--Horn theorem applied to the trace-zero matrix
\(A-\tr(A)y_0y_0^{\top}\), where \(y_0=A^{-1/2}x_0\).

For the total \((n-1)\)-dimensional measure \(S\) of the facets we prove
\[
       S_{\max}(\cE_A)=2^n n^{-(n-2)/2}\sqrt{\det A}\,\sqrt{\tr(A^{-1})}.
\]
For \(n\ge3\) the maximisers are more rigid: on the sphere they are orthotopes
with all edge lengths equal and with a Schur--Horn equal diagonal condition for
\(A^{-1}\).  The prescribed-vertex facet-area problem is therefore equivalent to
a restricted Schur--Horn problem with a prescribed barycentric basis.  In
dimension two this recovers the Connes--Zagier property for ellipses.  In
dimension three, however, the direct higher-dimensional analogue fails for
triaxial ellipsoids at principal-axis vertices; an exact obstruction is given.
\end{abstract}

\maketitle

\section{Introduction}

The study of extremal polytopes inscribed in convex bodies is classical.  A
particularly striking planar theorem, observed in a ballistic context by
Richard~\cite{Richard2004} and proved elegantly by Connes and Zagier
\cite{ConnesZagier}, says that an ellipse has the following property: through
every boundary point there passes an inscribed parallelogram of maximal
perimeter.  The maximal perimeter is \(4\sqrt{a^2+b^2}\), where \(a,b\) are the
semiaxes of the ellipse.

This paper studies the analogous questions for parallelepipeds in centred
ellipsoids in \(\R^n\).  We use the closed ellipsoid
\[
   \cE_A=\{x\in\R^n:x^{\top}A^{-1}x\le1\},
\]
where \(A\) is symmetric positive definite, and we call a full-dimensional
parallelepiped \emph{inscribed} if all its vertices lie on \(\partial\cE_A\).  Two
functionals are considered:
\begin{itemize}[leftmargin=2em]
\item the total length \(L\) of the one-skeleton;
\item the total \((n-1)\)-dimensional measure \(S\) of the \(2n\) facets.
\end{itemize}
The first functional is the perimeter when \(n=2\); the second functional is the
usual surface area when \(n=3\), and also coincides with the perimeter when
\(n=2\).

The basic parametrisation is simple but important.  If \(B=A^{1/2}\) denotes the
positive-definite square root of \(A\), then every inscribed parallelepiped is
centred at the origin and has edge vectors
\[
       v_i=\lambda_i B u_i\qquad (i=1,\ldots,n),
\]
where \(U=[u_1\ \cdots\ u_n]\in O(n)\), \(\lambda_i>0\), and
\(\sum_i\lambda_i^2=4\).  Thus, after applying \(B^{-1}\), the problem is one of
orthotopes inscribed in the unit sphere.

For the one-skeleton length, the Cauchy--Schwarz inequality gives the sharp
bound
\[
       L(P)\le 2^n\sqrt{\tr A}.
\]
The equality condition is flexible enough to solve the prescribed-vertex problem
completely.  If \(x_0\in\partial\cE_A\) and \(y_0=B^{-1}x_0\in S^{n-1}\), then an
orthonormal basis satisfying
\[
     \diag\bigl(U^{\top}(A-\tr(A)y_0y_0^{\top})U\bigr)=0
\]
produces an extremal parallelepiped through \(x_0\); such a basis exists by the
Schur--Horn theorem.  After changing signs of the columns of \(U\), all
coordinates of \(U^{\top}y_0\) are positive, and these coordinates determine the
edge lengths.  This gives a higher-dimensional edge-length analogue of the
Connes--Zagier phenomenon.

For the facet-area functional, the answer is different.  The global maximum is
\[
       S_{\max}(\cE_A)
       =2^n n^{-(n-2)/2}\sqrt{\det A}\,\sqrt{\tr(A^{-1})}.
\]
For \(n\ge3\), equality forces all spherical edge lengths to be equal and forces
\(\diag(U^{\top}A^{-1}U)\) to be the constant vector
\(\tr(A^{-1})\one/n\).  With a prescribed vertex, this means that the orthonormal
basis must also be barycentric with respect to \(y_0\).  Thus the constrained
facet-area problem becomes a restricted Schur--Horn problem on the stabiliser of
\(\one\).  In dimension two this again gives the Connes--Zagier theorem, because
facet area and perimeter coincide.  In dimension three the higher-dimensional
facet-area analogue is false for triaxial ellipsoids at principal-axis endpoints:
if \(y_0\) is an eigenvector of \(A^{-1}\), then an extremal surface-area
parallelepiped through \(x_0\) exists only when the two eigenvalues of
\(A^{-1}\) on \(y_0^{\perp}\) are equal.

Throughout, vectors are columns, \(\|\cdot\|\) is the Euclidean norm, and
\((\cdot)^{\top}\) denotes transpose.  For a square matrix \(M\),
\(\diag(M)\) denotes the vector of diagonal entries.  For a vector \(d\),
\(\Diag(d)\) denotes the diagonal matrix with diagonal \(d\).  If
\(z=(z_i)\), then \(z^{\circ 2}:=(z_i^2)\) is the componentwise square.

\section{Parametrising inscribed parallelepipeds}

We begin with the elementary reduction to the sphere.  It also explains why no
separate translation parameter appears in the later formulae.

\begin{lemma}[Centring and spherical parametrisation]\label{lem:parametrisation}
Let \(A\) be positive definite and let
\[
   P=c+\left\{\frac12\sum_{i=1}^n t_i v_i:-1\le t_i\le 1\right\}
\]
be a full-dimensional parallelepiped in \(\R^n\).  Suppose that all vertices of
\(P\) lie on \(\partial\cE_A\).  Then \(c=0\).  Moreover, with
\(B=A^{1/2}\), the edge vectors are of the form
\[
       v_i=\lambda_i B u_i\qquad (i=1,\ldots,n),
\]
where \(U=[u_1\ \cdots\ u_n]\in O(n)\), \(\lambda_i>0\), and
\(\sum_i\lambda_i^2=4\).  Conversely, every such choice gives an inscribed
parallelepiped in \(\cE_A\).
\end{lemma}

\begin{proof}
Put \(C=A^{-1}\), \(V=[v_1\ \cdots\ v_n]\), and let
\(\varepsilon=(\varepsilon_1,\ldots,\varepsilon_n)\in\{\pm1\}^n\).  The vertices
are
\[
       x_\varepsilon=c+\frac12 V\varepsilon,
\]
and the hypothesis is
\[
       x_\varepsilon^{\top}C x_\varepsilon=1
       \qquad(\varepsilon\in\{\pm1\}^n).
\]
After expanding, the left-hand side is
\[
 c^{\top}Cc+\varepsilon^{\top}V^{\top}Cc+\frac14\varepsilon^{\top}V^{\top}CV\varepsilon.
\]
This polynomial in the signs is constant on the discrete cube.  Hence all its
non-constant Walsh coefficients vanish.  Equivalently, multiplying by
\(\varepsilon_i\), respectively by \(\varepsilon_i\varepsilon_j\), and averaging
over all sign vectors gives
\[
       (V^{\top}Cc)_i=0,
       \qquad
       (V^{\top}CV)_{ij}=0\quad(i\ne j).
\]
Since \(V\) is invertible and \(C\) is positive definite, \(V^{\top}Cc=0\)
implies \(c=0\).  Thus the edge vectors are pairwise orthogonal for the inner
product \(\ip{x}{y}_C=x^{\top}Cy\), and the remaining constant condition is
\[
       \frac14\sum_{i=1}^n v_i^{\top}Cv_i=1.
\]
Now write \(w_i=B^{-1}v_i\).  Then
\(\ip{w_i}{w_j}=v_i^{\top}Cv_j\), so the non-zero vectors \(w_i\) are Euclidean
orthogonal and satisfy \(\sum_i\|w_i\|^2=4\).  Setting
\(\lambda_i=\|w_i\|\) and \(u_i=w_i/\lambda_i\) gives the stated
parametrisation.  The converse is immediate: for every sign vector,
\[
       \left\|\frac12\sum_{i=1}^n \varepsilon_i\lambda_i u_i\right\|^2
       =\frac14\sum_{i=1}^n\lambda_i^2=1,
\]
so all vertices of the spherical orthotope lie on \(S^{n-1}\), and applying
\(B\) puts all vertices on \(\partial\cE_A\).
\end{proof}

We shall use the following standard form of the Schur--Horn theorem.

\begin{theorem}[Schur--Horn]\label{thm:schur-horn}
Let \(M\) be a real symmetric \(n\times n\) matrix with eigenvalue vector
\(\mu\), written in non-increasing order.  A vector \(d\in\R^n\) occurs as
\(\diag(U^{\top}MU)\) for some \(U\in O(n)\) if and only if \(d\) is majorised by
\(\mu\).  In particular:
\begin{enumerate}[label=\normalfont(\alph*)]
\item for every symmetric \(M\), there is \(U\in O(n)\) such that
\[
       \diag(U^{\top}MU)=\frac{\tr M}{n}\one;
\]
\item if \(\tr M=0\), there is \(U\in O(n)\) such that
\[
       \diag(U^{\top}MU)=0.
\]
\end{enumerate}
\end{theorem}

\begin{proof}
The theorem is the classical Schur--Horn theorem; see, for example,
\cite[Theorem~4.3.26]{HornJohnson}.  The two stated consequences follow because
the constant vector with the same sum as the eigenvalues is majorised by the
eigenvalue vector.  In the trace-zero case this constant vector is \(0\).
\end{proof}

\section{The total length of the one-skeleton}\label{sec:length}

Let \(P\) be an inscribed parallelepiped with edge data
\(v_i=\lambda_iBu_i\) as in Lemma~\ref{lem:parametrisation}.  There are
\(2^{n-1}\) edges parallel to \(v_i\).  Hence the total length of the one-skeleton
is
\begin{equation}\label{eq:L-formula}
       L(P)=2^{n-1}\sum_{i=1}^n\|v_i\|
       =2^{n-1}\sum_{i=1}^n\lambda_i\sqrt{u_i^{\top}Au_i}.
\end{equation}

\begin{theorem}[Maximal total edge length]\label{thm:L-global}
For every parallelepiped \(P\) inscribed in \(\cE_A\),
\[
       L(P)\le L_{\max}(\cE_A):=2^n\sqrt{\tr A}.
\]
Equality holds if and only if, for the corresponding orthonormal basis
\(U=[u_1\ \cdots\ u_n]\),
\begin{equation}\label{eq:L-global-equality}
       \lambda_i=\frac{2\sqrt{u_i^{\top}Au_i}}{\sqrt{\tr A}}
       \qquad(i=1,\ldots,n).
\end{equation}
In particular, every orthonormal basis \(U\) gives a maximiser by the choice
\eqref{eq:L-global-equality}.
\end{theorem}

\begin{proof}
Set \(b_i=\sqrt{u_i^{\top}Au_i}\).  Since \(U\) is orthogonal,
\(\sum_i b_i^2=\tr A\).  By the Cauchy--Schwarz inequality and
\(\sum_i\lambda_i^2=4\),
\[
       \sum_{i=1}^n\lambda_i b_i
       \le \left(\sum_i\lambda_i^2\right)^{1/2}
            \left(\sum_i b_i^2\right)^{1/2}
       =2\sqrt{\tr A}.
\]
Multiplying by \(2^{n-1}\) gives the bound.  Equality in the
Cauchy--Schwarz inequality is equivalent to \(\lambda_i=tb_i\) for a constant
\(t>0\).  The condition \(\sum_i\lambda_i^2=4\) gives
\(t=2/\sqrt{\tr A}\).
\end{proof}

\subsection{The prescribed-vertex problem for total edge length}

Let \(x_0\in\partial\cE_A\), and put
\[
       y_0=A^{-1/2}x_0\in S^{n-1}.
\]
By changing signs of the edge directions, we may assume that the prescribed
vertex is the all-plus vertex.  Thus an inscribed parallelepiped with spherical
edge data \((U,\lambda)\) has vertex \(x_0\) precisely when
\begin{equation}\label{eq:vertex-condition}
       \sum_{i=1}^n\lambda_i u_i=2y_0.
\end{equation}
Equivalently, if \(z=U^{\top}y_0\), then
\begin{equation}\label{eq:lambda-from-z}
       z_i>0,
       \qquad
       \lambda_i=2z_i
       \qquad(i=1,\ldots,n).
\end{equation}

\begin{proposition}[Equality condition with a prescribed vertex]\label{prop:L-vertex-condition}
Let \(x_0\in\partial\cE_A\) and \(y_0=A^{-1/2}x_0\).  An inscribed
parallelepiped through \(x_0\), represented by \((U,\lambda)\) and satisfying
\eqref{eq:lambda-from-z}, attains \(L_{\max}(\cE_A)\) if and only if
\begin{equation}\label{eq:L-vertex-equality}
       \diag(U^{\top}AU)=\tr(A)\,z^{\circ2},
       \qquad z=U^{\top}y_0.
\end{equation}
Equivalently,
\begin{equation}\label{eq:trace-zero-L}
       \diag\bigl(U^{\top}(A-\tr(A)y_0y_0^{\top})U\bigr)=0.
\end{equation}
\end{proposition}

\begin{proof}
Using \eqref{eq:lambda-from-z} in \eqref{eq:L-formula}, we obtain
\[
       L(P)=2^n\sum_{i=1}^n z_i\sqrt{u_i^{\top}Au_i}.
\]
Since \(z_i>0\) and \(\sum_i z_i^2=1\), the Cauchy--Schwarz inequality gives
\[
       L(P)\le 2^n
          \left(\sum_i z_i^2\right)^{1/2}
          \left(\sum_i u_i^{\top}Au_i\right)^{1/2}
       =2^n\sqrt{\tr A}.
\]
Equality holds exactly when the vectors
\((z_i)\) and \((\sqrt{u_i^{\top}Au_i})\) are proportional.  Squaring and using
\(\sum_i u_i^{\top}Au_i=\tr A\) gives \eqref{eq:L-vertex-equality}.  The
equivalence with \eqref{eq:trace-zero-L} follows from
\(z_i=\ip{u_i}{y_0}\).
\end{proof}

\begin{theorem}[Prescribed vertices for total edge length]\label{thm:L-prescribed}
For every \(x_0\in\partial\cE_A\), there exists an inscribed parallelepiped
\(P\subset\cE_A\) having \(x_0\) as a vertex and satisfying
\[
       L(P)=L_{\max}(\cE_A)=2^n\sqrt{\tr A}.
\]
\end{theorem}

\begin{proof}
Let \(T=\tr A\), \(y_0=A^{-1/2}x_0\), and
\[
       H=A-Ty_0y_0^{\top}.
\]
Then \(\tr H=0\).  By Theorem~\ref{thm:schur-horn}, there is
\(U=[u_1\ \cdots\ u_n]\in O(n)\) such that \(\diag(U^{\top}HU)=0\).  Put
\(z=U^{\top}y_0\).  No coordinate of \(z\) is zero: if \(z_i=0\), then
\[
       0=u_i^{\top}Hu_i=u_i^{\top}Au_i,
\]
contradicting positive definiteness of \(A\).  Replacing \(u_i\) by
\(\operatorname{sign}(z_i)u_i\) if necessary, we may assume that all coordinates
of \(z\) are positive, while the identity \(\diag(U^{\top}HU)=0\) is unchanged.
Now set \(\lambda_i=2z_i\).  Then \eqref{eq:vertex-condition} holds, so the
resulting parallelepiped has vertex \(x_0\), and
Proposition~\ref{prop:L-vertex-condition} gives equality in the length bound.
\end{proof}

\begin{corollary}[The planar Connes--Zagier property]\label{cor:CZ}
When \(n=2\), for every \(x_0\in\partial\cE_A\) there exists an inscribed
parallelogram through \(x_0\) with maximal perimeter.  The maximal perimeter is
\[
       4\sqrt{\tr A}.
\]
\end{corollary}

\begin{proof}
In dimension two, the total length of the one-skeleton is the perimeter.  The
statement is therefore Theorem~\ref{thm:L-prescribed} with \(n=2\).
\end{proof}

\section{The total facet area}\label{sec:facet}

Let again \(v_i=\lambda_iBu_i\), and set
\[
       V=[v_1\ \cdots\ v_n],
       \qquad
       \Lambda=\Diag(\lambda_1,\ldots,\lambda_n),
       \qquad
       C=A^{-1}.
\]
The Gram matrix of the edge vectors is
\[
       G=V^{\top}V=\Lambda U^{\top}AU\Lambda.
\]
For \(i=1,\ldots,n\), let \(G_{\widehat{i}}\) denote the principal submatrix of
\(G\) obtained by deleting row and column \(i\).  The two facets parallel to all
edge directions except \(v_i\) have \((n-1)\)-dimensional volume
\(\sqrt{\det G_{\widehat{i}}}\).  Hence
\begin{equation}\label{eq:S-gram}
       S(P)=2\sum_{i=1}^n\sqrt{\det G_{\widehat{i}}}.
\end{equation}
Since \(G\) is positive definite,
\[
       \det G_{\widehat{i}}=(\det G)(G^{-1})_{ii}.
\]
Moreover,
\[
       \det G=(\det\Lambda)^2\det A,
       \qquad
       G^{-1}=\Lambda^{-1}U^{\top}CU\Lambda^{-1}.
\]
Substituting these identities into \eqref{eq:S-gram} gives
\begin{equation}\label{eq:S-key}
       S(P)=2\sqrt{\det A}\left(\prod_{j=1}^n\lambda_j\right)
       \sum_{i=1}^n
       \frac{\sqrt{(U^{\top}CU)_{ii}}}{\lambda_i}.
\end{equation}

We shall need one elementary inequality for the edge-length parameters.

\begin{lemma}\label{lem:lambda-phi}
Let \(\lambda_i>0\) and \(\sum_i\lambda_i^2=4\).  Then
\begin{equation}\label{eq:Phi-bound}
       \left(\prod_{i=1}^n\lambda_i\right)
       \left(\sum_{i=1}^n\frac1{\lambda_i^2}\right)^{1/2}
       \le 2^{n-1}n^{(2-n)/2}.
\end{equation}
If \(n\ge3\), equality holds if and only if
\(\lambda_1=\cdots=\lambda_n=2/\sqrt n\).  If \(n=2\), equality holds for every
positive pair \((\lambda_1,\lambda_2)\) with \(\lambda_1^2+\lambda_2^2=4\).
\end{lemma}

\begin{proof}
Put \(t_i=\lambda_i^2/4\).  Then \(t_i>0\) and \(\sum_i t_i=1\).  Squaring the
left-hand side of \eqref{eq:Phi-bound} gives
\[
       4^{n-1}\left(\prod_{i=1}^n t_i\right)
       \left(\sum_{i=1}^n\frac1{t_i}\right)
       =4^{n-1}e_{n-1}(t_1,\ldots,t_n),
\]
where \(e_{n-1}\) is the elementary symmetric polynomial of degree \(n-1\).  By
Maclaurin's inequality,
\[
       e_{n-1}(t_1,\ldots,t_n)
       \le \binom{n}{n-1}\left(\frac{t_1+\cdots+t_n}{n}\right)^{n-1}
       =n^{2-n}.
\]
This proves the bound.  For \(n\ge3\), the equality case in Maclaurin's
inequality forces \(t_1=\cdots=t_n=1/n\).  For \(n=2\), the quantity
\(e_1(t_1,t_2)=t_1+t_2\) is identically equal to \(1\).
\end{proof}

\begin{theorem}[Maximal total facet area]\label{thm:S-global}
For every parallelepiped \(P\) inscribed in \(\cE_A\),
\begin{equation}\label{eq:Smax}
       S(P)\le S_{\max}(\cE_A)
       :=2^n n^{-(n-2)/2}\sqrt{\det A}\,\sqrt{\tr(A^{-1})}.
\end{equation}

If \(n\ge3\), equality holds if and only if
\begin{equation}\label{eq:S-equality-nge3}
       \lambda_1=\cdots=\lambda_n=\frac2{\sqrt n}
       \quad\text{and}\quad
       \diag(U^{\top}A^{-1}U)=\frac{\tr(A^{-1})}{n}\one.
\end{equation}
Such maximisers exist.

If \(n=2\), the functional \(S\) is the perimeter.  In this case equality in
\eqref{eq:Smax} is equivalent to equality in Theorem~\ref{thm:L-global}; the
value is
\[
       S_{\max}(\cE_A)=4\sqrt{\det A}\sqrt{\tr(A^{-1})}=4\sqrt{\tr A}.
\]
\end{theorem}

\begin{proof}
Let
\[
       a_i=\sqrt{(U^{\top}CU)_{ii}}.
\]
Since \(\sum_i a_i^2=\tr C\), the Cauchy--Schwarz inequality applied to
\eqref{eq:S-key} gives
\[
       S(P)
       \le 2\sqrt{\det A}\sqrt{\tr C}
       \left(\prod_{j=1}^n\lambda_j\right)
       \left(\sum_{i=1}^n\frac1{\lambda_i^2}\right)^{1/2}.
\]
Lemma~\ref{lem:lambda-phi} gives \eqref{eq:Smax}.

Assume first that \(n\ge3\).  Equality in Lemma~\ref{lem:lambda-phi} gives
\(\lambda_i=2/\sqrt n\) for all \(i\).  Equality in the Cauchy--Schwarz
inequality then gives that all \(a_i\) are equal.  Since
\(\sum_i a_i^2=\tr C\), this is exactly
\(\diag(U^{\top}CU)=\tr(C)\one/n\).  Conversely, these conditions plainly give
equality.  Existence follows from the Schur--Horn theorem,
Theorem~\ref{thm:schur-horn}, applied to \(C=A^{-1}\).

For \(n=2\), the total facet measure is the perimeter, hence equals \(L\).  The
maximisers are therefore those described in Theorem~\ref{thm:L-global}.  The
identity
\[
       \det(A)\tr(A^{-1})=\tr(A)
\]
for \(2\times2\) positive definite matrices gives the stated value.
\end{proof}

\begin{remark}\label{rem:n2-equality-different}
The equality statement for \(S\) is genuinely different in dimension two.  For
\(n\ge3\), maximal facet area forces equal spherical edge lengths.  For \(n=2\),
Lemma~\ref{lem:lambda-phi} is an equality for every admissible pair
\((\lambda_1,\lambda_2)\), and the maximisers are exactly the maximal-perimeter
parallelograms from Theorem~\ref{thm:L-global}.
\end{remark}

\section{Prescribed vertices for the facet-area functional}\label{sec:S-prescribed}

The constrained problem for \(S\) is now best understood by combining the global
equality statement with the vertex condition \eqref{eq:vertex-condition}.

\begin{definition}\label{def:barycentric}
Let \(y\in S^{n-1}\).  An orthonormal basis
\(U=[u_1\ \cdots\ u_n]\in O(n)\) is called \emph{barycentric with respect to}
\(y\) if
\begin{equation}\label{eq:barycentric}
       \ip{u_1}{y}=\cdots=\ip{u_n}{y}=\frac1{\sqrt n}.
\end{equation}
Equivalently,
\[
       \frac1{\sqrt n}\sum_{i=1}^n u_i=y.
\]
Such bases exist for every \(y\): choose any orthogonal matrix sending
\(n^{-1/2}\one\) to \(y\).
\end{definition}

\begin{proposition}[Prescribed-vertex equality for facet area]\label{prop:S-vertex-condition}
Let \(n\ge3\), \(x_0\in\partial\cE_A\), \(y_0=A^{-1/2}x_0\), and
\(C=A^{-1}\).  There exists an inscribed parallelepiped through \(x_0\) with
\[
       S(P)=S_{\max}(\cE_A)
\]
if and only if there exists an orthonormal basis \(U\) barycentric with respect
to \(y_0\) such that
\begin{equation}\label{eq:S-barycentric-condition}
       \diag(U^{\top}CU)=\frac{\tr C}{n}\one.
\end{equation}
When such a basis exists, the extremal parallelepiped is obtained by taking
\(\lambda_i=2/\sqrt n\) for all \(i\).
\end{proposition}

\begin{proof}
Suppose first that \(P\) is an inscribed parallelepiped through \(x_0\) and that
\(S(P)=S_{\max}(\cE_A)\).  By Theorem~\ref{thm:S-global}, all
\(\lambda_i=2/\sqrt n\) and \eqref{eq:S-barycentric-condition} holds.  After
changing signs of edge directions, the prescribed vertex is the all-plus vertex,
so \eqref{eq:vertex-condition} gives
\[
       y_0=\frac1{\sqrt n}\sum_{i=1}^n u_i.
\]
Thus \(U\) is barycentric with respect to \(y_0\).

Conversely, if such a barycentric \(U\) is given and \(\lambda_i=2/\sqrt n\),
then \eqref{eq:vertex-condition} holds, so the parallelepiped has vertex
\(x_0\).  The equality condition in Theorem~\ref{thm:S-global} gives
\(S(P)=S_{\max}(\cE_A)\).
\end{proof}

\begin{corollary}\label{cor:S-n2}
For \(n=2\) and every \(x_0\in\partial\cE_A\), there exists an inscribed
parallelogram through \(x_0\) satisfying
\[
       S(P)=S_{\max}(\cE_A)=4\sqrt{\tr A}.
\]
\end{corollary}

\begin{proof}
In dimension two, \(S\) is the perimeter.  The result follows from
Corollary~\ref{cor:CZ}.
\end{proof}

\begin{proposition}[Balls and ellipsoids of revolution]\label{prop:revolution}
Let \(n\ge3\), \(x_0\in\partial\cE_A\), and \(y_0=A^{-1/2}x_0\).
If
\begin{equation}\label{eq:revolution}
       A=a\,y_0y_0^{\top}+b\,(\Id-y_0y_0^{\top})
       \qquad(a,b>0),
\end{equation}
then there is an inscribed parallelepiped through \(x_0\) that simultaneously
attains \(L_{\max}(\cE_A)\) and \(S_{\max}(\cE_A)\).  In particular, this holds for
Euclidean balls for every prescribed boundary point.
\end{proposition}

\begin{proof}
Choose a basis \(U=[u_1\ \cdots\ u_n]\) barycentric with respect to \(y_0\), and
set \(\lambda_i=2/\sqrt n\).  Then \(P\) has vertex \(x_0\).  For each \(i\),
\(\ip{u_i}{y_0}^2=1/n\), and therefore
\[
       u_i^{\top}Au_i=\frac{a}{n}+b\left(1-\frac1n\right)=\frac{\tr A}{n}.
\]
Thus \(\diag(U^{\top}AU)=\tr(A)\one/n=\tr(A)(U^{\top}y_0)^{\circ2}\), so the
length equality condition \eqref{eq:L-vertex-equality} holds.  The inverse
matrix has the same form,
\[
       A^{-1}=a^{-1}y_0y_0^{\top}+b^{-1}(\Id-y_0y_0^{\top}),
\]
so similarly
\(\diag(U^{\top}A^{-1}U)=\tr(A^{-1})\one/n\).  Proposition~\ref{prop:S-vertex-condition}
then gives equality for \(S\) as well.
\end{proof}

The next result shows that the facet-area analogue of the Connes--Zagier
property does not hold for general ellipsoids in dimension three.

\begin{theorem}[A three-dimensional obstruction]\label{thm:3d-obstruction}
Let \(n=3\), let \(C=A^{-1}\), and suppose that
\(y_0\in S^2\) is an eigenvector of \(C\).  Let \(\gamma_2,\gamma_3\) be the two
eigenvalues of \(C\) on \(y_0^{\perp}\).  There exists an inscribed
parallelepiped through \(x_0=A^{1/2}y_0\) attaining \(S_{\max}(\cE_A)\) if and only
if
\[
       \gamma_2=\gamma_3.
\]
Consequently, a triaxial ellipsoid in \(\R^3\) has no maximal-surface-area
inscribed parallelepiped through an endpoint of a principal axis.
\end{theorem}

\begin{proof}
If \(\gamma_2=\gamma_3\), then \(A\) has the form \eqref{eq:revolution} with
axis \(y_0\), and Proposition~\ref{prop:revolution} applies.

Conversely, assume that an \(S\)-maximiser through \(x_0\) exists.  By
Proposition~\ref{prop:S-vertex-condition}, there is a basis
\(U=[u_1,u_2,u_3]\) barycentric with respect to \(y_0\) and satisfying
\(\diag(U^{\top}CU)=\tr(C)\one/3\).  Choose orthonormal coordinates
\((y_0,e_2,e_3)\) in which
\[
       C=\Diag(\gamma_1,\gamma_2,\gamma_3).
\]
Since \(U\) is barycentric, each column has the form
\[
       u_i=\frac1{\sqrt3}y_0+h_i,
       \qquad h_i\in y_0^{\perp},
\]
where the three vectors \(h_i\) have norm \(\sqrt{2/3}\) and pairwise inner
product \(-1/3\).  Thus they are the vertices of a regular triangle in
\(y_0^{\perp}\).  For a suitable angle \(\theta\), we may write
\[
       h_i=\sqrt{\frac23}
       \bigl(\cos(\theta+2\pi(i-1)/3)e_2+
             \sin(\theta+2\pi(i-1)/3)e_3\bigr).
\]
The equality of the three diagonal entries of \(U^{\top}CU\) gives equality of
\[
       \gamma_2\cos^2(\theta+2\pi(i-1)/3)
       +\gamma_3\sin^2(\theta+2\pi(i-1)/3),
       \qquad i=1,2,3.
\]
If \(\gamma_2\ne\gamma_3\), this would force the three numbers
\(\cos^2(\theta+2\pi(i-1)/3)\) to be equal.  Equivalently, the three numbers
\(\cos(2\theta+4\pi(i-1)/3)\) would be equal.  Their sum is zero, so they would
all have to be zero, which is impossible for three angles separated by
\(2\pi/3\).  Hence \(\gamma_2=\gamma_3\).
\end{proof}

\section{The restricted Schur--Horn problem with barycentre}\label{sec:restricted}

Proposition~\ref{prop:S-vertex-condition} leaves a natural finite-dimensional
problem.  Fix a positive definite matrix \(C\) and a unit vector \(y\).  Define
\[
   \mathcal B(C,y)=
   \left\{U\in O(n):
      U \text{ is barycentric with respect to } y,
      \ \diag(U^{\top}CU)=\frac{\tr C}{n}\one
   \right\}.
\]
For \(n\ge3\), the prescribed-vertex maximum for the facet-area functional is
attained at \(x_0=A^{1/2}y\) if and only if \(\mathcal B(A^{-1},y)\ne\varnothing\).

This is a restricted version of the Schur--Horn theorem.  Indeed, fix one
barycentric basis \(U_0\) for \(y\).  Then all barycentric bases for \(y\) are of
the form \(U_0V\), where
\[
       V\in O(n),\qquad V\one=\one.
\]
Thus the question is whether, for
\(M=U_0^{\top}CU_0\), one can find an orthogonal matrix \(V\) in the stabiliser
of \(\one\) such that
\[
       \diag(V^{\top}MV)=\frac{\tr M}{n}\one.
\]
The unrestricted Schur--Horn theorem guarantees such a diagonal after conjugacy
by an arbitrary orthogonal matrix; the difficulty is the additional constraint
\(V\one=\one\).  Theorem~\ref{thm:3d-obstruction} shows that this additional
constraint is substantial.

\begin{problem}\label{prob:characterise-B}
For \(n\ge3\), characterise the pairs \((C,y)\), with \(C\) positive definite
and \(y\in S^{n-1}\), for which \(\mathcal B(C,y)\ne\varnothing\).
Equivalently, characterise the boundary points of an ellipsoid through which
there passes an inscribed parallelepiped of globally maximal facet area.
\end{problem}

The geometric prescribed-vertex problem is complete in the following cases: in
dimension two, it is always soluble because it is the Connes--Zagier perimeter
problem; in all dimensions, it is soluble for ellipsoids of revolution at a
point on the axis of revolution; and in dimension three at a principal-axis
point it is soluble exactly in the rotationally symmetric case described in
Theorem~\ref{thm:3d-obstruction}.

\subsection*{Competing interests} There are no competing interests to disclose.

\subsection*{Acknowledgements}
The author is grateful to Dominik Bysiewicz and Marcin Guzik for inspiring
discussions on the topic. 

\subsection*{Funding}
Funding received from the Institute of Mathematics of Czech Academy of Sciences (RVO 67985840) is acknowledged with thanks.


\begin{thebibliography}{9}

\bibitem{ConnesZagier}
A.~Connes and D.~Zagier,
A property of parallelograms inscribed in ellipses,
\emph{Amer. Math. Monthly} \textbf{114} (2007), no.~10, 909--914.

\bibitem{HornJohnson}
R.~A. Horn and C.~R. Johnson,
\emph{Matrix Analysis}, 2nd ed.,
Cambridge Univ. Press, 2013.

\bibitem{Richard2004}
J.-M.~Richard,
Safe domain and elementary geometry,
\emph{Eur. J. Phys.} \textbf{25} (2004), 835--844.

\end{thebibliography}
\end{document}